\documentclass[12pt]{article}

\usepackage[utf8]{inputenc}
\usepackage[english]{babel}
\usepackage{amsmath}
\usepackage{amssymb}
\usepackage{amscd}
\usepackage{amsthm}
\usepackage{dsfont}
\usepackage{mathtools}
\usepackage{graphicx}
\usepackage{epstopdf}
\usepackage{hyperref}
\usepackage{cite}
\usepackage{authblk}
\usepackage[top=2cm, bottom=2cm, left=2cm, right=2cm]{geometry}
\usepackage{indentfirst}
\usepackage{xcolor}
\usepackage{euscript}
\mathtoolsset{
showonlyrefs=true,
mathic=true,
}

\allowdisplaybreaks

\hypersetup{backref,
 colorlinks=false}
\hypersetup{pdfborder=0 0 0}

\newcommand{\erfc}[1] {\mathrm{Erfc}\ #1 }

\renewcommand{\epsilon}{\varepsilon}
\renewcommand{\kappa}{\varkappa}
\renewcommand{\phi}{\varphi}

\renewcommand{\Im}[1]{\mathrm{Im}\ #1}

\newcommand{\sign}[1] {\,\mathrm{sign} \!\left(#1\right) }

\newcommand{\mexp}[2][]{\mathbb E^{#1}\! \left[#2 \right]}

\newcommand{\ltransinv}[2][\!\!]{\mathfrak{L}^{-1}_{#1} \left[ #2 \right]}

\newtheorem*{mythm}{Theorem}

\begin{document}
\title{Skew Brownian motion with dry friction: The~Pugachev--Sveshnikov approach}

\author{Sergey Berezin}
\author{Oleg Zayats}

\affil{Department of Applied Mathematics\\ Peter the Great St.Petersburg Polytechnic University, Russia \\ \href{mailto:servberezin@yandex.ru}{servberezin@yandex.ru}, \href{mailto:zay.oleg@gmail.com}{zay.oleg@gmail.com}}

\date{}
\maketitle 

\begin{abstract}
The Caughey--Dieness process, also known as the Brownian motion with two valued drift, is used in theoretical physics as an advanced model of the Brownian particle velocity if the resistant force is assumed to be dry friction. This process also appears in many other fields, such as applied physics, mechanics, astrophysics, and pure mathematics. In the present paper we are concerned with a more general process, skew Brownian motion with dry friction. The probability distribution of the process itself and of its occupation time on the positive half line are studied. The approach based on the Pugachev--Sveshnikov equation is used.
\end{abstract}

\section{Introduction}
\label{sec:intro}
Brownian motion plays an important role in statistical physics and other applied areas of science. The classical physical theory of Brownian motion was developed by Einstein and Smoluchovwski\cite{Einstein} in the beginning of 20th century.
In this theory, a weightless Brownian particle is driven by the random force resulting from the collision with molecules. Also, the friction between this particle and the molecules, in other words the resistant force due to collision, is assumed to be viscous, that is, proportional to the particle's velocity. Mathematically, these assumptions lead to the Wiener process. The major drawback of this model is that almost all the trajectories of the particle are continuous but nowhere differentiable, so the velocity cannot be defined properly. Later on, this problem was solved by Ornstein and Uhlenbeck~\cite{Ornstein}, who developed a refined physical theory that accounts for the particle's inertia.

A more evolved model of Brownian motion uses an assumption that along with viscous, dry friction is present, which is independent of the particle's speed, but depends on the motion direction. Such a model in its simplest form, when the viscous friction is absent, is called Brownian motion with dry friction and  was first studied by Caughey and Dienes~\cite{Caughey} in 1960s. They considered a process similar to Ornstein--Uhlenbeck's but the linear term was replaced by the sign function. This process turned out to be useful in analysis of different phenomena in control theory~\cite{Fuller}, seismic mechanics~\cite{Crandall}, communication systems theory~\cite{Lindsey}, radio physics~\cite{Tikhonov}, nonlinear stochastic dynamics~\cite{Kempe}, and also in pure mathematics~\cite{Shiryaev, Lejay2, Gairat}. In 2000s, other applications of the Caughey--Dienes process emerged. The major interest was in nanofrictional systems~\cite{Riedo2004}, particles separation~\cite{Eglin2006}, ratchets~\cite{Buguin2006, Fleishman2007}, granular motors~\cite{Talbot2011}, dynamics of granular media~\cite{Hayakawa2004, Hayakawa2005}, and in dynamics of droplets on moving surfaces~\cite{Gennes, Daniel2005, Mettu2010, Goohpattader2010}. Some additional publications on the subject can be found in the author's work~\cite{Berezin1}. 

The present paper deals with the model of skew Brownian motion with dry friction, or more precisely with the so-called skew Caughey--Dienes process, which can be understood as the velocity of a particle in this model. In what follows, for simplicity, we limit ourselves to studying skewing at zero only.

In order to understand how skewing works we consider an excursion of the process~$V(t)$; an excursion is the part of the process's trajectory, located between two time points at which~$V(t)$ vanishes, that is to say, between two consequent stops of a particle. An excursion of the Caughey--Dienes process without skewing is known to be a sufficiently good approximation of the real Brownian particle behavior. However, note that for such a process, due to symmetry, positive and negative excursions appear with the same probability~$0.5$, which is not always the case. Generally speaking, in reality the probability~$\alpha$ to have a positive excursion is different from~$0.5$, that is, some skewing takes place. In this way we come to the model of the skew Caughey--Dieness process, where the probability~$\alpha$ is an additional parameter of the model, which has to be estimated from the experimental data. Note that due to the above interpretation of the skewing, it is not that difficult to do.

Skew processes are tightly connected with diffusion in discontinuous media and has numerous applications in several fields of physics~\cite{Lejay}. For example, they describe shock acceleration of charged particles in a magnetic field~\cite{Zhang2000} and diffusion in geophysical problems associated with the study of inhomogeneous porous media~\cite{Ramirez2006}. It should also be noted that such processes arise in some special randomly perturbed Hamiltonian systems~\cite{Freidlin1994}. Besides, the skewing procedure was studied from purely mathematical standpoint for a number of typical stochastic processes. The most detailed results are available for the process of skew Brownian motion~\cite{Lejay}, for the Ornstein--Uhlenbeck process~\cite{Wang2015}, and for the Cox--Ingersoll--Ross process~\cite{Tian2018}. Unfortunately, the skew Caughey--Dienes process is not that well known, and this, in particular, motivates our study.

Generally, calculating the probabilistic characteristics of a skew diffusion is considered to be a difficult mathematical problem~\cite{Lejay2}. The standard tools to solve it are the Fokker--Planck--Kolmogorov equation, Feynman-Kac equation, or the random walk approximation. We propose an alternative way, based on the Pugachev--Sveshnikov equation, which was used in authors' previous works~\cite{Berezin1, Berezin2}. From our perspective, this approach allows one to get to the result faster, and in a more algorithmic way.

\section{Skew Caughey--Dienes process}
\label{sec:base-section}
Skewing is intimately connected with the notion of local time, and it is possible to proof that the skew Caughey--Dienes process described in the introduction can be represented as a unique strong solution~\cite{Lejay} of the following stochastic differential equation
\begin{equation}
  \label{eq:sde1}
  dX(t) = -2 \mu \sign{X(t)} dt + (2 \alpha - 1) dL_X^0(t) + \sqrt{2} dW(t),\, t>0,\quad X(0)=0,
\end{equation}
where~$\alpha \in (0,1)$ is the skewing parameter, the probability of an excursion to be positive. By~$W(t)$ we denote a standard Wiener process starting at zero, and $L_X^0(t)$ is the symmetric local time of~$X(t)$ at the level zero:
\begin{equation}
  L_X^0(t) = \lim_{\epsilon \to +0} \frac{1}{2 \epsilon} \int_0^t \mathds{1}_{(-\epsilon, +\epsilon)}(X(s))\, d[X]_s,
\end{equation}
where~$[X]_s =2 s$ is the so-called quadratic variation of~$X(s)$. Further, we will be interested in the positive half-line occupation time of~$X(t)$
\begin{equation}
  \label{eq:occup_time}
  \mathcal{I}(t) = \int_0^t \mathds{1}_{(0, +\infty)}(X(s))\, ds.
\end{equation}

One can think of~$X(t)$ and~$\mathcal{I}(t)$ as of the components of the vector diffusion process~$(X(t),\mathcal{I}(t))$ governed by the system of SDEs
\begin{equation}
  \label{eq:syst1}
  \left\{
  \begin{aligned}
    &dX(t) = -2 \mu \sign{X(t)} dt +(2 \alpha - 1) dL_X^0(t) + \sqrt{2} dW(t),\\
    &d\mathcal{I}(t) = \mathds{1}_{(0, +\infty)}(X(s))\, ds,\quad X(0)=\mathcal{I}(0)=0.
  \end{aligned}
\right.
\end{equation}

In what follows, we derive explicit formulas for the probability density function of~\eqref{eq:syst1}, following ideas from~\cite{Berezin2}. 

It can be shown~\cite{Berezin2} that the characteristic function~$E(z_1,z_2;t)=\mexp{e^{i (z_1 X(t) + z_2 \mathcal{I}(t)) }}$ of the process~$(X(t), \mathcal{I}(t))$ satisfies the equation
\begin{equation}
  \label{eq:ps_eq}
  \frac{\partial E}{\partial t} + (z_1^2 - i z_2 /2) E + (2 \mu z_1 - z_2/2) \hat{E} - 2 i (2 \alpha - 1) z_1 \Psi_0= 0, \quad E(z_1,z_2;0) = 1,
\end{equation}
where we adopt the short notation~$\hat{E}(z_1,z_2;t)$ and~$\Psi_0(z_2,t)$:
\begin{equation}
  \label{eq:notations}
  \hat{E} = \frac{1}{\pi} \mathrm{v.p.} \int\limits_{- \infty}^{+\infty} \frac{E|_{z_1=s}}{s-z_1}\, ds, \quad \Psi_0 = \frac{1}{2 \pi} \mathrm{v.p.} \int\limits_{- \infty}^{+\infty} E|_{z_1=s}\, ds.
\end{equation}

For~$\Im{\zeta}\ne 0$ let us introduce the Cauchy-type integral~$\Phi(\zeta,z_2; t)$ and its limit values~$\Phi^{\pm}(z, z_2;t)$  on the real axis from upper and lower half-planes (with respect to the first argument):
\begin{equation}
\label{eq:Cauchy_type_int}
    \Phi(\zeta,z_2;t)=\frac{1}{2 \pi i}\int \limits_{-\infty}^{+\infty} \frac{E|_{z_1=s}}{s-\zeta} ds,\ \Phi^{\pm} (z,z_2; t) = \lim_{\zeta \to z\pm i0} \Phi(\zeta, z_2; t), \ \Im{z} = 0.
\end{equation}
It is well known that~$\Phi(\cdot, z_2; t)$ is analytic when~$\Im{\zeta}\ne 0$, and that~$\Phi^\pm$ satisfy the Sokhotski--Plemelj formulas when~$\Im{z} = 0$:
\begin{equation}
  \label{eq:sokh}
  \Phi^+ - \Phi^- = E, \ \Phi^+ + \Phi^- = - i \hat{E}.
\end{equation}
Clearly, one can rewrite~\eqref{eq:ps_eq} in terms of~$\Phi^\pm$, that gives the condition of the Riemann boundary value problem, and now we need to recover analytic functions~$\Phi^\pm$, from this condition. Applying the Laplace transform with respect to~$t$, and denoting its argument by~$p$, we get to the formula
\begin{equation}
  \label{eq:rbvp}
  (z_1^2 + 2\mu i z_1 + p - i z_2) \tilde{\Phi}^+ -i (2 \alpha - 1) z_1 \tilde{\Psi}_0- \frac{1}{2} = (z_1^2 - 2\mu i z_1 + p) \tilde{\Phi}^- + i (2 \alpha - 1) z_1 \tilde{\Psi}_0 + \frac{1}{2}.
\end{equation}
The Laplace transforms are labeled with tildes above the functions.

 Note that the left-hand side of~\eqref{eq:rbvp} can be analytically continued for all~$z_1 \in \mathbb{C}$ such that~$\Im{z_1}>0$, also the right-hand side can be analytically continued for all~$z_1 \in \mathbb{C}$ such that~$\Im{z_1}<0$. Since they match when~$\Im{z_1} = 0$, they turn out to be the elements of the same entire function of argument~$z_1 \in \mathbb{C}$ . Assuming that~$\Phi^{\pm} (z, z_2; t) =  O(\frac{1}{|z|})$ when~$z \to \infty$ for~$\Im{z} \gtrless 0$, by generalized Liouville's theorem one can realize that this entire function is actually linear: $G_0(z_2, t) + z_1 G_1(z_2, t)$. This leads to the equality
\begin{equation}
  \label{eq:sol_chf}
  \tilde{\Phi}^\pm = \frac{G_0 + z_1 G_1 \pm i (2 \alpha - 1) z_1 \tilde{\Psi}_0 \pm 1/2}{z_1^2 \pm 2 i \mu z_1 +p  - (1 \pm 1) i z_2/2}.
\end{equation}

Note that the denominator in~\eqref{eq:sol_chf} has zeros~$i \nu^{\pm} = i(- \mu \pm \sqrt{\mu^2 + p - i z_2})$ and~$i \kappa^{\pm} = i( \mu \pm \sqrt{\mu^2 + p})$ such that~$\Im{\!\!(i \nu^{\pm})} \gtrless 0$ and $\Im{\!\!(i \kappa^{\pm})} \gtrless 0$. At the same time, $\tilde{\Phi}^\pm$ should be analytic in upper and lower half-planes, therefore, the singularities at~$i \nu^+$ and~$i \kappa^-$ are removable. This gives a system of linear equations to determine~$G_0$ and~$G_1$. Also, taking into account the definition of~$\Psi_0$ in~\eqref{eq:notations} and performing integration in~\eqref{eq:sol_chf} one can find that~$\tilde{\Psi}_0 = - i G_1$. After that, the system of linear equations for~$G_0$ and~$G_1$ can be written in the following form
\begin{equation}
  \label{eq:g01}
  G_0 + 2 i \nu^+ \alpha G_1 + 1/2 =0,\ G_0 + 2 i \kappa^- (1- \alpha)G_1 - 1/2=0.
\end{equation}

 Finally, substituting~$G_0$ and~$G_1$ from~\eqref{eq:g01} into~\eqref{eq:sol_chf}, one can get to the final expression for~$\tilde{E}$, using the first of the Sokhotski--Plemelj formulas~\eqref{eq:sokh}. Particularly, after necessary simplifications one can get the Laplace transform of the characteristic function of~$X(t)$ and~$\mathcal{I}(t)$:
 \begin{equation}
   \label{im_ch}
  \begin{aligned}
    &\tilde{E}_X(z_1,p) = \tilde{E}(z_1,0; p) =  \frac{1}{i \kappa^-} \left( \frac{\alpha}{z_1 + i \kappa^+} - \frac{1- \alpha}{z_1 - i \kappa^+} \right),\\
    &\tilde{E}_\mathcal{I}(z_2,p) = \tilde{E}(0,z_2; p) =  \frac{\alpha \kappa^+-(1-\alpha) \nu^-}{\nu^- \kappa^+((1-\alpha) \kappa^- - \alpha \nu^+)}.
  \end{aligned}
\end{equation}

Straightforward computations using the table of Laplace and Fourier transforms for the first expression in~\eqref{im_ch} gives us the probability density function of~$X(t)$:
\begin{equation}
 f_X(x, t) = \left(\frac{1}{\sqrt{\pi t}} e^{-\frac{(|x|+2 \mu t)^2}{4t}} + \mu e^{-2 \mu |x|} \erfc{\frac{|x|-2 \mu t}{2 \sqrt{t}}} \right) \cdot \left\{
   \begin{aligned}
     &\alpha, &&x>0,\\
     &1-\alpha, &&x<0.
   \end{aligned}
 \right.
\end{equation}

Letting~$t$ tend to~$+\infty$, one can obtain the steady-state probability density function
\begin{equation}
 f_X^{\infty}(x) = f_X(x, +\infty) =2 \mu e^{-2 \mu |x|}\left\{
   \begin{aligned}
     &\alpha, &&x>0,\\
     &1-\alpha, &&x<0.
   \end{aligned}
 \right.
\end{equation}

\begin{figure}[h!]
\centering
\includegraphics[width=0.5\linewidth]{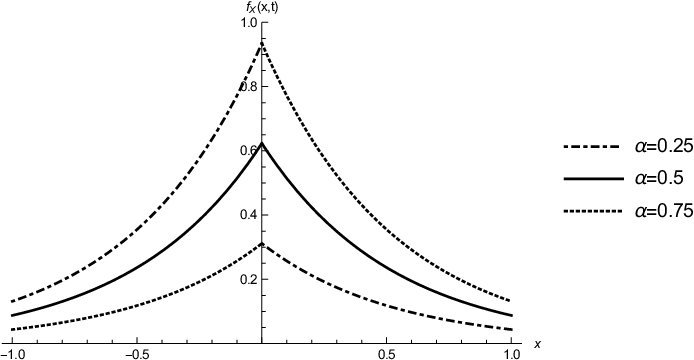}
\caption{Probability density function~$f_X(x, t)$ for different~$\alpha$ ($\mu=1, t=1$).}
\label{fig1}
\end{figure}

The plot of the probability density function of~$X(t)$ for different~$\alpha$ is given in Fig.~\ref{fig1}. Note that for~$\alpha \ne 0.5$ the curve has a jump at zero, so that the probability that~$X(t)>0$ is equal to~$\alpha$.

Handling the second line in~\eqref{im_ch} takes somewhat more effort. First, we introduce the function:
\begin{equation}
  \label{s_rational}
  S(w_1,w_2) = \frac{\alpha (w_2 + \mu) + (1-\alpha) (w_1 + \mu)}{(w_1+\mu)(w_2+\mu)(\alpha (w_1 - \mu) + (1-\alpha) (w_2-\mu))},
\end{equation}
then the expression for~$\tilde{E}_\mathcal{I}(z_2,p)$ can be written as 
\begin{equation}
  \tilde{E}_\mathcal{I}(z_2,p) = S(\sqrt{\mu^2 + p - i z_2}, \sqrt{\mu^2 + p}).
\end{equation}

After that, we note that since~$\mathcal{I}$ is non-negative, instead of using the inverse Fourier transform with respect to~$z_2$, one can use the inverse Laplace transform with respect to~$q = - i z_2$:
\begin{equation}
  f_{\mathcal{I}}(y, t) = \ltransinv[p,t]{\ltransinv[q,y]{S(\sqrt{ \mu ^2 + p + q},\sqrt{\mu^2 + p})}},
\end{equation}
where~$\ltransinv[p,t]{\cdot}$ is the inverse Laplace transform with respect to the parameter~$p$, and the argument of the original is~$t$.

Now we write the following chain of equalities:
\begin{align}
  \label{eq:cd_occup_time_eq20_1}
  &f_{\mathcal{I}}(y, t)= \ltransinv[p,t]{\ltransinv[q,y]{S(\sqrt{ \mu ^2 + p+q},\sqrt{\mu^2 + p}) }}=e^{ - \mu^2 t} \ltransinv[p,t]{\ltransinv[q,y]{S(\sqrt{p+q}, \sqrt{p}) }}\\
  \label{eq:cd_occup_time_eq20_3}
  &=e^{ -\mu^2 t} \ltransinv[p,t]{e^{-p y} \ltransinv[q,y]{S(\sqrt{q}, \sqrt{p})}} =\chi(t-y) e^{ - \mu^2 t} \ltransinv[p, t-y]{ \ltransinv[q,y]{S(\sqrt{q}, \sqrt{p})}}\\
  \label{eq:cd_occup_time_eq20_5}
  &=\frac{\chi(t-y) e^{- \mu^2 t}} {4 \pi (y(t-y))^{3/2}} \int \limits_0^{+ \infty} \int \limits_0^{+ \infty} \ltransinv[p,s_2] { \ltransinv[q,s_1]{S(q,p)}} \ s_1 s_2\ e^{ -\frac{s_1^2}{4y}-\frac{s_2^2}{4(t-y)}}\ ds_1 ds_2\\
  \label{eq:cd_occup_time_eq20_8}
  &=\frac{4\chi(t-y) e^{- \mu^2 t}} {\pi \sqrt{y(t-y)}} \int \limits_0^{+ \infty} \int \limits_{0}^{+ \infty} \ltransinv[p,2\sqrt{t-y}s_2]{\ltransinv[q,2 \sqrt{y} s_1]{S(q,p)}}\ s_1 s_2\ e^{-s_1^2 -s_2^2}\ ds_1 ds_2.
\end{align}

In the first and second lines we used well-known properties of the Laplace transform, in the line three we applied Efros's theorem~\cite{Shabat}, and in the last line we made a change of variables in the double integral.

As the function~$S$ in~\eqref{s_rational} is rational, it is easy to calculate~$\ltransinv[p,2\sqrt{t-y}s_1]{\ltransinv[q,2 \sqrt{y} s_2]{S(p,q)}}$. Skipping cumbersom but trivial in nature computations, we get to the final formula for the probability density function of~$\mathcal{I}$
\begin{equation}
  f_{\mathcal{I}}(y, t) = \frac{4 e^{- \mu^2 t}} {\pi \sqrt{y(t-y)}} \int \limits_0^{+ \infty} \int \limits_{0}^{+ \infty}\, \chi(2 \sqrt{y}s_1, 2 \sqrt{t-y}s_2)\, s_1 s_2 e^{-s_1^2 - s_2^2}\, ds_1 ds_2, \quad 0<y<t,
\end{equation}
where~$\erfc{\!(\cdot)}$ is the complementary error function, and
\begin{equation}
  \begin{aligned}    
    &\chi(s_1,s_2) = \frac{1-\alpha}{\alpha} e^{-\mu \left(s_1 - \frac{2-\alpha}{\alpha}s_2\right)} \chi^+(s_1,s_2) + \frac{\alpha}{1-\alpha} e^{\mu \left(\frac{1+\alpha}{1-\alpha}s_1 - s_2\right)} \chi^-(s_1,s_2),\\
    &\chi^+(s_1,s_2) = \mathds{1}_{(0,+\infty)}(\alpha s_1 - (1-\alpha)s_2), \quad   \chi^-(s_1,s_2) = 1-\chi^+(s_1,s_2).
  \end{aligned}
\end{equation}

\begin{figure}[h!]
\centering
\includegraphics[width=0.5\linewidth]{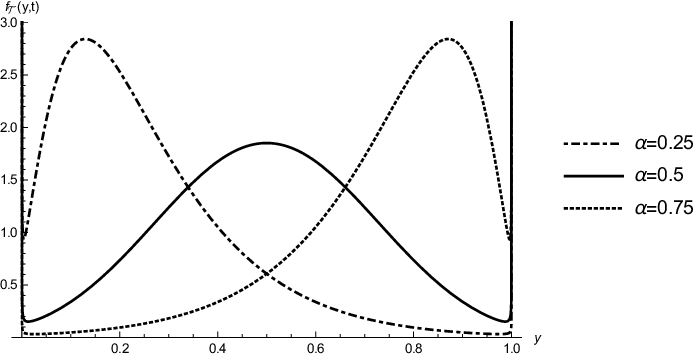}
\caption{Scaled occupation time density function~$f_{\mathcal{T}}(y, t)$ for different~$\alpha$ ($\mu=1, t=2$).}
\label{fig2}
\end{figure}

Not performing any other simplifications of this formula, we illustrate the generic shape of the scaled occupation time~$\mathcal{T}(t) = \mathcal{I}(t)/t$ in Fig.~\ref{fig2}. The probability density function of~$\mathcal{T}(t)$ is given by~$f_{\mathcal{T}}(y, t) = t f_{\mathcal{I}}(t y, t)$.

It is worth noticing that unlike the well-known arcsine law, the distribution in Fig.~\ref{fig2} is unimodal, and its mode can be controlled by the parameter~$\alpha$. This completely agrees with the physical interpretation of the skew Caughey--Dienes process given in the introduction.  

At the end of the section, we formulate our result in the form of the following theorem.
\begin{mythm}
  The PDF of~$X(t)$, the steady-state PDF of~$X(t)$, and the PDF of the positive half-line occupation time~$\mathcal{I}(t)$ are given by:
  \begin{equation}
 f_X(x, t) = \left(\frac{1}{\sqrt{\pi t}} e^{-\frac{(|x|+2 \mu t)^2}{4t}} + \mu e^{-2 \mu |x|} \erfc{\frac{|x|-2 \mu t}{2 \sqrt{t}}} \right) \cdot \left\{
   \begin{aligned}
     &\alpha, &&x>0,\\
     &1-\alpha, &&x<0,
   \end{aligned}
 \right.
\end{equation}
\begin{equation}
 f_X^{\infty}(x) = f_X(x, +\infty) =2 \mu e^{-2 \mu |x|} (\alpha \mathds{1}_{(0, +\infty)}(x) +  (1-\alpha) \mathds{1}_{(-\infty,0]}(x)),
\end{equation}
\begin{equation}
  f_{\mathcal{I}}(y, t) = \frac{4 e^{- \mu^2 t}} {\pi \sqrt{y(t-y)}} \int \limits_0^{+ \infty} \int \limits_{0}^{+ \infty}\, \chi(2 \sqrt{y}s_1, 2 \sqrt{t-y}s_2)\, s_1 s_2 e^{-s_1^2 - s_2^2}\, ds_1 ds_2, \quad 0<y<t
\end{equation}
where~$\erfc{\!(\cdot)}$ is the complementary error function, and
\begin{equation}
  \begin{aligned}    
    &\chi(s_1,s_2) = \frac{1-\alpha}{\alpha} e^{-\mu \left(s_1 - \frac{2-\alpha}{\alpha}s_2\right)} \chi^+(s_1,s_2) + \frac{\alpha}{1-\alpha} e^{\mu \left(\frac{1+\alpha}{1-\alpha}s_1 - s_2\right)} \chi^-(s_1,s_2),\\
    &\chi^+(s_1,s_2) = \mathds{1}_{(0,+\infty)}(\alpha s_1 - (1-\alpha)s_2), \quad   \chi^-(s_1,s_2) = 1-\chi^+(s_1,s_2).
  \end{aligned}
\end{equation}
\end{mythm}

\section{Conclusions}
We derived explicit formulas for the probability density function of the skew Caughey--Dienes process and its occupation time on the positive half-line, which generalizes known results for the regular Caughey--Dienes process. In fact, more general result was obtained for the Laplace transform of the joined characteristic function. Essentially, our approach is based on the reduction to a Riemann boundary value problem, and clearly it can be used to find the characteristics of more general SDEs with piecewise linear coefficients and local time.


\begin{thebibliography}{9}

\bibitem{Einstein} \emph{Einstein~A., Smoluchowski~M.} Brownian motion. Collected papers.~--- M.-L.: ONTI.~--- 1936 (in Russian).

\bibitem{Ornstein} \emph{Uhlenbeck~G.~E., Ornstein~L.~S.} On the theory of Brownian motion // Physical Review.~--- 1930.~--- Vol.~36, no.~5.~--- P.~823--841.

\bibitem{Caughey} \emph{Caughey~T.~K., Dienes~J.~K.} Analysis of a nonlinear first-order system with a white noise input // Journal of Applied Physics.~--- 1961.~--- Vol.~32, no.~11.~--- P.~2476--2479.

\bibitem{Fuller} \emph{Fuller~A.~T.} Exact analysis of a first-order relay control system with a white noise disturbance // International Journal of Control.~--- 1980.~--- Vol.~31, no.~5.~--- P.~841--867.

\bibitem{Crandall} \emph{Crandall~S.~H., Lee~S.~S., Williams~J.~H.} Accumulated slip of a friction-controlled mass excited by earthquake motions // Journal of Applied Mechanics.~--- 1974.~--- Vol.~41, no.~4.~---P.~1094--1098.

\bibitem{Lindsey} \emph{Lindsey~W.C.} Synchronization systems in communication and control.~--- N.Y.: Prentice Hall.~--- 1972.

\bibitem{Tikhonov} \emph{Tikhonov~V.~I., Mironov~M.~A.} Markov processes.~--- M.: Sovetskoe radio. ~--- 1977 (in Russian).

\bibitem{Kempe} \emph{Ahlbehrendt~N., Kempe~V.} Analysis of stochastic systems. Nonlinear dynamical systems.~--- Berlin: Akademie Verlag.~--- 1984 (in German).

\bibitem{Gairat} \emph{Gairat~A., Scherbakov~V.} Density of skew Brownian motion and its functionals with applications in finance // Mathematical Finance.~--- 2016.~--- Vol.~27, no.~4.~--- P.~1069--1088.
  
\bibitem{Lejay2} \emph{Lejay~A., Len\^otre~L., Pichot~G.} One-dimensional skew diffusions: explicit expressions of densities and resolvent kernels // [Research Report] Inria Rennes - Bretagne Atlantique; Inria Nancy - Grand Est. 2015. ~--- 2015. P.~ 1--28.
  
\bibitem{Shiryaev} \emph{Shiryaev~A.~N., Cherny~A.~S.} Some distributional properties of a Brownian motion with a drift and an extension of P.~L\'evy's theorem // Theory of Probability \& Its Applications.~--- 2000.~--- Vol.~44, no.~2.~--- P.~412--418.

\bibitem{Riedo2004} \emph{Riedo~E., Gnecco~E.} Thermally activated effects in nanofriction // Nanotechnology.~--- 2004.~--- Vol.~15, no.~4.~--- P.~S288--S292.

\bibitem{Eglin2006} \emph{Eglin~M., Eriksson~M. A., Carpick~R. W.} Microparticle manipulation using inertial forces // Applied Physics Letters.~--- 2006.~--- Vol.~88, no.~9.~--- P.~2172401.

\bibitem{Buguin2006} \emph{Buguin~A., Brochard~F., de Gennes~P.-G.} Motions induced by asymmetric vibrations // The European Physical Journal E.~--- 2006.~--- Vol.~19, no.~1.~--- P.~31--36.

\bibitem{Fleishman2007} \emph{Fleishman~D., Asscher~Y., Urbakh~M.} Directed transport induced by asymmetric surface vibrations: making use of friction // Journal of Physics: Condensed Matter.~--- 2007.~--- Vol.~19, no.~9.~--- P.~096004.

\bibitem{Talbot2011} \emph{Talbot~J., Wildman~R. D., Viot~P.} Kinetics of a frictional granular motor // Physical Review Letters.~--- 2011.~--- Vol.~107, no. 13.~--- P.~138001.
\bibitem{Hayakawa2005} \emph{Hayakawa~H.} Langevin equation with Coulomb friction // Physica D.~--- 2005.~--- Vol.~205, no.~1.~--- P.~48--56.

  
\bibitem{Hayakawa2004} \emph{Kawarada~A., Hayakawa~H.} Non-Gaussian velocity distribution function in a vibrating granular bed // Journal of the Physical Society of Japan.~--- 2004.~--- Vol.~73, no.~8.~--- P.~2037--2040.
  
\bibitem{Daniel2005} \emph{Daniel~S., Chaudhury~M. K., de Gennes~P.-G} Vibration-actuated drop motion on surfaces for batch microfluidic processes // Langmuir.~--- 2005.~--- Vol.~21, no.~9.~--- P.~4240--4248.

\bibitem{Mettu2010} \emph{Mettu~S., Chaudhury~M. K.} Stochastic relaxation of the contact line of a water drop on a solid substrate subjected to white noise vibration: Roles of hysteresis // Langmuir.~--- 2010.~--- Vol.~26, no.~11.~--- P.~8131--8140.

\bibitem{Goohpattader2010} \emph{Goohpattader~P. S., Chaudhury~M. K.} Diffusive motion with nonlinear friction: apparently Brownian // The Journal of Chemical Physics.~--- 2010.~--- Vol.~133, no.~2.~--- P.~3460530.
  
\bibitem{Gennes} \emph{de Gennes~P.~-G.} Brownian motion with dry friction // Journal of Statistical Physics.~--- 2005.~--- Vol.~119, no.~5--6.~--- P.~953--962.



\bibitem{Berezin1} \emph{Berezin~S., Zayats~O.} Energy dissipation in a friction-controlled slide of a body excited by random motions of the foundation // Physical Review E.~---2018.~--- Vol.~97, iss.~1.--- P.~012144.
  
\bibitem{Lejay} \emph{Lejay~A.} On the constructions of the skew Brownian motion // Probability Surveys.~--- 2006.~--- no.~3.~--- P.~413--466.

\bibitem{Zhang2000} \emph{Zhang~M.} Calculation of diffusive shock acceleration of charged particles by skew Brownian motion // Astrophysical Journal.~--- 2000.~--- Vol.~541.~--- P.~428--435.

\bibitem{Ramirez2006} \emph{Ramirez~J. M., Thomann~E. A., Waymir~E.C., Haggerty~R., Wood~B.} A generalization of Taylor--Aris  formula and skew diffusion // Multiscale Modeling and Simulation.~--- 2006.~--- Vol.~5, no.~3.~--- P.~786--801.


\bibitem{Freidlin1994} \emph{Freidlin~M. I., Wentzel~A. D.} Random perturbation of Hamiltonian systems // Memoirs of the American Mathematical Society.~--- 1994.~--- Vol.~109, no.~523.~--- P.~1--87.

  
  
  
\bibitem{Wang2015} \emph{Wang~S., Song~S., Wang~Y.}  Skew Ornstein-Uhlenbeck processes and their financial applications // Journal of Computational and Applied Mathematics.~--- 2015.~--- Vol.~273, no.~1.~--- P.~363--382.

\bibitem{Tian2018} \emph{Tian~Y., Zhang~H.}  Skew CIR process, conditional characteristic function, moments and bond pricing // Applied Mathematics and Computation.~--- 2018.~--- Vol.~329.~--- P.~230--238.
  
\bibitem{Berezin2} \emph{Berezin~S., Zayats~O.} Application of the Pugachev--Sveshnikov equation to the Baxter occupation time problem // Informatics and Applications.~--- 2015.~--- Vol.~9, no.~2.~--- P.~39--49.

\bibitem{Shabat} \emph{Shabat~B.V.} Introduction to complex analysis.~--- M.: Nauka.~--- 1969 (in Russian).
  
\end{thebibliography}
\end{document}